\documentclass[11pt]{article}
\usepackage[utf8]{inputenc}
\usepackage[margin=1in,left=1in]{geometry}
\usepackage{setspace}
\usepackage{graphicx}
\usepackage{amssymb, amsmath, amsthm, graphics}
\usepackage{mathabx,epsfig}
\usepackage[mathscr]{euscript}
\usepackage{tikz}
\usepackage{pigpen}
\usepackage{bm}
\usepackage{mathtools}
\usetikzlibrary{calc}
\usetikzlibrary{matrix}
\usepackage{latexsym}
\usepackage[all]{xy}
\usepackage{color}

\usepackage{bbold}

\usepackage{tikz}
\input xy
\xyoption{all}
\usepackage{euscript}
\usepackage{multirow}
\usepackage{etex, pictexwd,dcpic}
\usetikzlibrary{positioning}
\usetikzlibrary{shapes.geometric}
\usetikzlibrary{shapes.misc}
\usetikzlibrary{calc}
\usetikzlibrary{positioning}

\usepackage{pgflibraryarrows}
\usepackage{pgflibrarysnakes}

\usetikzlibrary{trees} % LATEX and plain TEX 
\usetikzlibrary[trees] % ConTEXt

\newtheorem{theorem}{Theorem}
\newtheorem{definition}[theorem]{Definition}
\newtheorem{lemma}[theorem]{Lemma}

\newtheorem{proposition}[theorem]{Proposition}

\newtheorem{remark}{Remark}

\newtheorem{example}{Example}

\numberwithin{equation}{section}

\renewcommand{\(}{\begin{equation*}}
\renewcommand{\)}{\end{equation*}}
\newcommand{\bea}{\begin{eqnarray*}}
\newcommand{\eea}{\end{eqnarray*}}

\def\endofproof {\hfill{$\Box$}\\}

\def\L{\ensuremath{{\cal L}}}

\def\L{\ensuremath{{\cal L}}}

\newcommand{\beq}{\begin{equation}}
\newcommand{\eeq}{\end{equation}}

\newcommand{\into}{\hookrightarrow}

\newcommand{\theproof}{\noindent {\bf Proof.\ }}

\numberwithin{equation}{section}

\renewcommand{\(}{\begin{equation}}
\renewcommand{\)}{\end{equation}}

\def\1{{\bf 1}}

\def\<{\langle}
\def\>{\rangle}

\numberwithin{equation}{section}

\newcommand{\R}{\ensuremath{\mathbb R}}
\newcommand{\RR}{\ensuremath{\mathbb R}}

\newcommand{\ZZ}{\ensuremath{\mathbb Z}}
\newcommand{\Z}{\ensuremath{\mathbb Z}}

\newcommand{\BB}{\ensuremath{\mathbf B}}

\newcommand{\CC}{\ensuremath{\mathbb C}}

\newcommand{\chp}{\ensuremath{\mathscr{C}\mathrm{h}^{+}}}

\newcommand{\sset}{\ensuremath{s\mathscr{S}\mathrm{et}}}

\newcommand{\cartsp}{\mathscr{C}\mathrm{art}\mathscr{S}\mathrm{p}}

\newcommand{\map}{\mathrm{Map}}

\begin{document}

\title{Twisted smooth Deligne cohomology}

 \author{Daniel Grady and Hisham Sati
 \\
  }
%{hsati@pitt.edu}, 

\maketitle

\begin{abstract} 
%Deligne cohomology, being a differential cohomology, combines
%integral and de Rham cohomology. 
Deligne cohomology can be viewed as a differential refinement of 
integral cohomology, hence captures both topological and geometric information. 
On the other hand, it can be viewed as the simplest nontrivial version of 
a differential cohomology theory. While more involved differential cohomology
theories have been explicitly twisted, the same has not been done to 
Deligne cohomology, although existence is known at a general abstract level. 
We work out what it means to twist Deligne cohomology,  
by taking degree one twists of both integral cohomology and  
 de Rham cohomology. 
 We present the main properties of 
the new theory and illustrate its use with examples and applications.   
Given how versatile Deligne cohomology
has proven to be, we believe that this explicit and utilizable treatment of its 
twisted version will be useful. 
  \end{abstract}

 \tableofcontents
 
%%%%%%%%%%%%%%%% 
\section{Introduction}
%%%%%%%%%%%%%%%%

Twistings form an interesting feature of (generalized) cohomology theories.
While on general grounds existence   is established at the axiomatic/abstract level
via parametrized spectra (see \cite{MS} \cite{Units}), constructing such theories 
explicitly is often a nontrivial task (see \cite{ABG} \cite{SW} \cite{LSW} for recent 
 illustrations of this). By the same token, differential refinements of (generalized) 
 cohomology theories are also established at the general abstract level
(see \cite{BNV} \cite{BS} \cite{Urs}). However, again, explicit 
utilizable constructions require considerable work (see \cite{GS2} \cite{GS3}). 
Combining the above two aspects, one also has twisted differential cohomology
theories at the general yet somewhat abstract level \cite{BN}. The goal of this 
paper is to work out an explicit case of such a differential twisted theory. 

\medskip
We will be interested  mainly in twists of differential refinements of ordinary, i.e. integral, 
cohomology. This cohomology theory has smooth extension with various different realizations, 
including those of \cite{Del} \cite{CS} \cite{Ga} \cite{Bry} \cite{DL} \cite{HS} \cite{BKS} \cite{BB}. 
All of these realizations are in fact isomorphic (see \cite{SSu} \cite{BS}). 
We will use the smooth Deligne cohomology incarnation \cite{Del} \cite{Be} \cite{Gi} \cite{Ja} 
\cite{EV} \cite{Ga}. In terms of machinery, we will use the 
approach of simplicial presheaves and higher stacks 
\cite{Cech} \cite{Urs}\cite{FSS1} \cite{HQ} \cite{FSS2} . This has the virtue of being
amenable to generalizations  and allowing the use of powerful 
algebraic machinery.

\medskip
%Therefore, here we take a step back and 
%Therefore, we consider twists of Deligne cohomology which, in  
%a sense, is the simplest and first nontrivial differential cohomology theory. 
Twisted Deligne cohomology is understood to exist from general constructions 
 \cite{Bu} \cite{BN} and should be in some sense
one of the simplest twisted differential cohomology theories. However,
we have not seen this theory discussed in detail anywhere in the literature.
This, together with the versatility and utility of Deligne cohomology, motivated us to 
believe that this would be a very useful task. Furthermore, we view this,  in some sense, as the
toy model and the prototype example for more involved twisted differential spectra
described this way. From the fact that constructing this theory was not as straightforward 
as one might have thought, with unexpected structures and subtleties, the 
involved task of unraveling the details turned out to be worthwhile.

%\medskip
%What would it mean to twist Deligne cohomology? At the practical level, Deligne cohomology is 
%a fusion of de Rham cohomology and {\v C}ech cohomology. One could naively 
%consider twists of de Rham cohomology, made periodic, as well as twists of {\v C{ech cohomology 
%separately. 
%%and then envision combining together for form a twisted Deligne operator 
%%$D + {\hat H}$, where $\hat{H}$ somehow co
%%mbines the twists from both de Rham
%%and {\v C}ech. 
%%\medskip
%The operator in Deligne cohomology of the {\v C}ech-de Rham double complex is 
%a combination of the de Rham differential $d$, acting in the direction of the de Rham
%complex, and the {\v C}ech differential 
%$\delta$, acting in the direction of the {\v C}ech complex. This the Deligne or {\v C}ech-de Rham operator
%$D=(-1)^p d + \delta$ (see \cite{BT} for an extensive discussion). 
%One could then consider two situations (or stages) to twisting the  operator $D$:
%First, adding a closed differential form $H$, i.e. modifying only the form part, leading to twisting of the
%der Rham differential $d \leadsto d_H$. Second, adding a differential cohomology class $\hat{H}$, i.e. 
%modifying both parts of the whole {\v C}ech-de Rham differential, i.e. $d \leadsto d_H$ and 
%$\delta \leadsto \delta_{B}$ where $(H, B)$ are appropriate (de Rham, {\v C}ech)-components of $\hat{H}$.
%We will see that it will not quite work that way, but this turns out to be a good heuristic to keep 
%in mind.  

\medskip
The study of the twisted de Rham cohomology is essential for understanding the 
confluent hypergeometric integral which can be regarded as a pairing of the twisted 
de Rham cohomology and a twisted cycle.
For the case of 1-forms there has been a considerable amount of work in this direction,
e.g. Deligne, who gives the twisted de Rham theorem in \cite[II; 6.3]{De} \cite{DM} as well as
\cite{Ki1} \cite{Ki2} \cite{AKOT} \cite{ASp}. Deligne's work has also had other incarnations, such 
as in Witten's approach to Morse theory \cite{Wi-Morse}.  
%We consider the twists for Deligne cohomology. Since Deligne cohomology is 
%a combination of de Rham cohomology (which calculates the real cohomology 
%of the manifold $M$) and integral cohomology, it makes sense that the twists will 
%be a combination of the integral twists and twists via differential forms. More 
%precisely, we will seek a differential form which models the rationalization 
%of an integral twist. 
Together with twisted integral cohomology, discussed
in Sec. \ref{Sec Z}, 1-form twisted de Rham cohomology lead to compatible
 twistings of Deligne cohomology in Sec. \ref{Sec twd}. Appropriately, the theory turns out to have
 pleasant computability properties that, again, are there axiomatically and abstractly, but that we 
 do explicitly verify and cast in utilizable form  in Sec. \ref{Sec Prop}. 
 We give several examples in Sec. \ref{Sec Ex} to illustrate the constructions.

\medskip
We now highlight  directly-related constructions in the literature. 
In \cite{Bu} differential refinements of integral cohomology are considered, leading to a sheaf-theoretic 
definition of smooth Deligne cohomology. 
In \cite{BKS} a bordism model for
the differential extension of ordinary integral cohomology is given in which one has integration
and products and a simple verification the projection formula.
Twistings of integral cohomology, as explained in \cite{Fr}, 
a priori are 1-dimensional, determined by a local system
$Z \to M$, which is a bundle of groups isomorphic to $\Z$. This is then determined
up to isomorphism by an element of $H^1(M; {\rm Aut}(\Z))\cong H^1(M; \Z/2)$.
The twisted integral cohomology $H^\bullet (M; \Z)$ may be defined using a 
cochain complex. This also admits a {\v C}ech description. We will make use of 
a generalization of this degree 1 local system setup in Sec. \ref{Sec twd}. 
An extensive discussion of the degree 1 case can be found in \cite{Farb}.

\medskip
Note that the Deligne complex can be viewed from more than one angle. From an algebraic 
point of view, it is the resolution of the group $\Z$. From a geometric (and topological) 
point of view, the starting point is the de Rham complex and this is viewed as the 
resolution of some sheaf. In the simplest case, this is $\R$, while more  involved situations
corresponding to twists, will require more delicate local systems.
%as we will see in Sec. \ref{Sec two} and Sec. \ref{Sec Main}. 
Note also  that there is a related concept of $\Z/2$-twisted de Rham forms \cite{HZ}, 
whereby  one can also twist by $\Z/2$ classes related to orientation; ${\cal O}_X$ is taken 
as the principal $\Z/2$-bundle of orientations of 
$TX$, and form ${\cal O}$-twisted $k$-forms as sections of 
${\cal O}\otimes_{\Z/2} \Lambda^k T^*X$. This is used to study twisted currents 
associated to Stiefel-Whitney classes. 

\medskip
We will be in the setting of homotopy sheaves \cite{Jardine} (see also \cite{Du} for a very 
readable account). In Sec. \ref{Sec twd}, we show that pulling back the universal bundle over 
a map  which classifies a twist, we get a bundle ${\cal H}^q\to M$ over $M$. We define the 
$\omega$-twisted Deligne cohomology of $M$ of degree $q$ to be the the connected components
$\pi_0\Gamma(M,{\cal H}^q)$.  
%This is particularly important when  we consider the higher degree twists in \cite{td}. 
We will mostly be dealing  with the category of sheaves of chain complexes,
while occasionally considering 
twisted differential cohomology within smooth sheaves of spectra.
%seems to have pleasant conceptual advantages. 
%, equipped with connection. 
%\medskip
%In  Sec. \ref{Sec Del} twisted Deligne cohomology groups of a smooth manifold $M$ 
%are defined using smooth stacks. This approach is also more suited to the higher twists 
%and is easily adapted to other differential function spectra, such as differential $K$-theory. 
%Moreover, this approach is easily adapted to other underlying theories which do 
%not have infinite loop spaces which are strictly abelian. 

 \medskip
%The approach we take will be based on chain complexes, as we indicated above. 
Chain complexes provide a useful way to present $H\ZZ$-module
spectra \cite{Shi}. 
Let ${\rm Ch}$ be the ordinary symmetric monoidal category of chain complexes 
of vector spaces. Formally invert the class of quasi-isomorphisms in Ch gives an $\infty$-category 
${\rm Ch}_\infty:=N({\rm Ch})[W^{-1}]$, the stable $\infty$-category obtained by 
localization of the $\infty$-category at quasi-isomorphisms.
The natural map $\iota: {\rm Ch} \to {\rm Ch}_\infty$ is a lax 
symmetric monodical functor. Furthermore, there is an
equivalence of symmetric monoidal categories \cite{Lu}  
$$
H: {\rm Ch}_\infty \buildrel{\sim}\over{\longrightarrow} {\rm Mod}_{H\Z}\;,
$$
where ${\rm Mod}_{H\Z}$ denotes the module spectra over the Eilenberg-MacLane spectrum 
$H\Z$. 
Consider the group of integers $\Z$ as an object in  $N({\bf Ch})[W^{-1}]$
by viewing it as a chain complex concentrated in degree zero. 
Being a stable $\infty$-category, ${\rm Ch}_\infty$ is
enriched over spectra. Hence the Eilenberg-MacLane spectrum can be
defined by the mapping spectrum 
$H\Z:= {\rm Map}(\Z, \Z)$.
Moreover, $H$ can be chosen such that $H(\Z)=H\Z$. 
As explained in \cite{Bu}, this can be considered as a commutative algebra 
in $N({\rm Sp})[W^{-1}]^{}$ (inverting stable 
equivalences on spectra) so that we can form its module category ${\rm Mod}_{H\Z}$.
The homotopy groups of $H\Z$ are given by
$$
\pi_*(H\Z)\cong 
\left\{
\begin{array}{cc}
\Z, & *=0\\
0, & *\neq 0.
\end{array}
\right.
$$
% \cite{Bu} (Fact 4.37) a 
 A differential refinement of a commutative ring spectrum $R$ is a triple
$(R, A, c)$ consisting of a CDGA $A$ over $\R$ together with an equivalence 
$$
c: R \wedge H\R \buildrel{\sim}\over{\longrightarrow} HA
$$ 
in ${\rm CAlg}({\rm Mod}_{H\R})$.
Shipley has shown that one can model every $H\R$-algebra by a CDGA \cite{Shi}. 
When $R$ has a differential refinement $\widehat{R}$ whose 
underlying CDGA is the graded ring $\pi_*(R) \otimes \R$ with trivial differential, 
$R \wedge H \R$ is called formal. In this case, there is an equivalence $c$ which is uniquely 
determined up to homotopy by the property that it induces the canonical identification on 
homotopy groups \cite{Bu}. 
This occurs for $H\Z$, for which one can choose a real model whose underlying CDGA is $\R$ 
concentrated in degree 0 \cite{BN}. 
 Ordinary cohomology is an example of what Bunke and Nikolaus \cite{BN} call 
a differentially simple spectrum. These have the property that there is a very good choice of 
a differential extension as well as control of differential twists. 
Assume that $R$ is a differentially simple spectrum and $(R, A, c)$ is the 
canonical differential extension with $A=\pi_*(R) \otimes \R$ and equivalence 
$c: R \wedge H\R \buildrel{\sim}\over{\longrightarrow} HA$ in 
${\rm CAlg}({\rm Mod}_{H\R})$. Then every topological $R$-twist $E$ 
on a manifold $M$ admits a differential refinement 
%$(E, {\cal M}, d)$, 
which is unique up to canonical equivalence   \cite[Theorem 9.5]{BN}.

\medskip
%In order to have a more complete picture of the twists of Deligne cohomology that 
%we constructed,  we provide {\v C}ech cocycle data in sec. \ref{Sec Cech}. The language 
%of smooth stacks provide a convenient language to obtain cocycle data from smooth 
%analogues of classifying spaces \cite{Cech}. In Sec. \ref{Sec Cech Z} we obtain 
%{\v C}ech cocycle data involving the higher twists for periodic integral cohomology 
%from Sec.\ref{Sec Per Z}.  
For ordinary integral cohomology, the twists are classified by $B\ZZ/2$. The pullback of a map $\eta:M\to B\ZZ/2$ by the universal $\ZZ$-bundle over $B\ZZ/2$ gives a $\ZZ$-bundle over $M$. The space of sections of this bundle is an infinite loop space and represents the $\eta$-twisted cohomology. Alternatively, we can think of the pullback bundle as a locally constant sheaf on $M$. From this point of view, the constant stack $\BB\ZZ/2\simeq \underline{B\ZZ/2}$ classifies all locally trivial $\ZZ$-bundles over $M$. The local sections of the bundle in this case form a \emph{sheaf of spectra} over $M$. The cohomology represented by this spectrum can be calculated explicitly by replacing a smooth manifold $M$ with its 
{\v C}ech nerve $C(\{U_{\alpha}\})$. We can identify the group of connected components 
of the mapping space $\map(C(\{U_{\alpha}\}),\BB \ZZ/2)$ with the {\v C}ech cohomology group 
$H^1(M;\ZZ/2)$
as the group  classifying the twists. To obtain cocycle data for a twist, we unravel
the bundle data coming from the action of $\ZZ/2$ on $\ZZ$. 
%Similarly, in Sec. \ref{Sec Cech Del} we obtain {\v C}ech cocycles for the general 
%twisted periodic Deligne cohomology that we constructed in Sec. \ref{Sec Higher}. 
%
%\medskip
%Summarizing, in the untwisted case, we have the following objects
%\begin{center}
%\begin{tabular}{|c||c|}
%\hline
% {\bf Untwisted cohomology} & {\bf Ordinary} 
%\\
%\hline\hline
%Underlying theory & Locally constant sheaf $\underline{\RR}$ 
%\\
%\hline
%de Rham complex & Ordinary de Rham complex $\Omega^*$ 
%\\
%\hline
%\end{tabular}
%\end{center}
%
%\medskip
%\noindent We wish to replace these with the twisted analogues
%\begin{center}
%\begin{tabular}{|c||c|}
%\hline
% {\bf Twisted cohomology} & {\bf Ordinary} 
% \\
% \hline
% \hline
% Underlying theory  & Locally constant sheaf ${\cal L}^{\eta}$ \\
%\hline
%Twist degree & One 
%\\
%\hline
%Twisting geometric object & Line bundle w/ flat conn. $d+H_1$ 
%\\
%\hline
%de Rham complex &  $\Omega^*(-; {\cal L}^{\eta} )$: differential $d+H_1$ \\
%\hline
%\end{tabular}
%\end{center}

%\medskip
%One might wonder if higher degree twists are possible. Indeed, 
%we will develop higher twisted
%Deligne cohomology in \cite{td}. This requires introducing a 
%periodic version of the theory.
%We will consider the higher twists in \cite{td}. 

\medskip
It has been brought to the authors' attention that a very interesting variation on the 
concept of twisted Deligne cohomology has been considered in the algebraic setting. 
Deligne-Beilinson cohomology with coefficients in a unipotent variations of 
mixed Hodge structure (VMHS) are considered in \cite{CH}. In \cite{H2} Hain, 
motivated by Hodge theory and motives,  develops Deligne-Beilinson cohomology 
of affine groups with a mixed Hodge structure. Kapranov has described MHS's 
in more geometric terms via certain categories of bundles with connections \cite{Kap}. 
Our paper should be considered as an approach to twisting Deligne cohomology
 from the point of view of differential geometry and algebraic topology. 

\medskip
The construction of surface holonomy of a bundle gerbe 
on unoriented surfaces and orientifolds 
can be described using twists for gerbes
\cite{SSW} \cite{GSW}  \cite{Gom} \cite{HMSV}. 
%These can be viewed as implicit constructions in 
 %twisted Deligne cohomology in the sense that  
 Translating our constructions to the language of 
 gerbes by writing out the {\v C}ech double complex corresponding 
 to our sheaves should recover the cocycle data discussed in  
the above works.  More  directly using our language of stacks is the 
 much more general model for ``higher orientifolds" given in 
 \cite[Sec. 4.4]{FSS3} to describe involutions arising from 
 M-theory on a manifold with boundary. In the same way that 
 Deligne cohomology is equivalent to gerbes with connections
 upon unraveling simplicial and cocycle data, our description of
 twisted Deligne cohomology should 
 likewise  be equivalent to twisted gerbes with connection in 
 all the above works.  The main point is that these two points 
 of view give  different  models for twisted differential integral cohomology.

\medskip
\par{\large \bf Acknowledgement.}
 The authors would like to thank the organizers and participants of the Geometric Analysis and 
 Topology Seminar at the Courant Institute for Mathematical Sciences for asking about twisting 
 Deligne cohomology, during a talk by H.S., which encouraged the authors to revisit and 
 carry out this project. D.G. would like to thank the Mathematics Department at the University of 
 Pittsburgh for hospitality during the final writing of this paper. The authors thank Richard Hain 
 for bringing to their attention related works in algebraic geometry
 and the referee for useful remarks. 

%%%%%%%%%
\section{Twists of integral cohomology}
\label{Sec Z}
%%%%%%%%%%

Integral cohomology appears as part of the data of ordinary differential cohomology. Therefore, 
twisting the former might give us some insight into the latter. Note, however, that finding the 
right way to do the twist is highly nontrivial, as was demonstrated in \cite{BN}. In this section, 
we review several models for twisted integral cohomology. The final approach uses the machinery 
of smooth stacks and fits in to the general machinery for twisted differential cohomology 
established in \cite{BN}. This last approach will be used in the 
discussion for twisted Deligne cohomology in subsequent sections.

\medskip
We consider twists of integral cohomology at the level of the 
Eilenberg-MacLane 
spectrum $H\Z$ as follows \cite{MQRT} \cite{ABG}. 
The corresponding infinite loop space is 
$\Omega^\infty H\Z\simeq K(\Z, 0) \simeq \Z$. This implies that the group of units is
${\rm GL}_1(H\Z)\simeq \{\pm 1\}\simeq \Z/2$, the invertible elements in $\Z$. 
Delooping then gives $B{\rm GL}_1(H\Z)\simeq B \Z/2 \simeq K(\Z/2, 1)$. 
Now consider a space $X$ with a twist given by a map
$X \to B{\rm GL}_1(H\Z)$. This can be viewed as an obstruction 
to orientation of a vector bundle $E$ over $X$ with respect to 
singular cohomology. This is given by the composite map that 
factors via the $J$-homomorphism through $B{\rm GL}_1(\mathbb{S})$. Here
$\mathbb{S}$ is the sphere spectrum, which is the unit for any spectrum,
including $H\Z$, so that there is always a map $\mathbb{S} \to H\Z$. That 
composite obstruction map is 
$$
\xymatrix{
X \ar[r]^E & BO \ar[r]^-{BJ} & B{\rm GL}_1(\mathbb{S}) \ar[r] &
B{\rm GL}_1(H\Z) \simeq K(\Z/2, 1)\;.
}
$$ 
This class $X \to K(\Z/2, 1)$ is, in fact, the first Stiefel-Whitney class
$w_1(X)$.

\medskip
The above homotopic description of the twist can be described geometrically 
as follows \cite{Fr}. One-dimensional twists of integral cohomology  is given by 
a local system $Z \to M$, which is  a bundle of groups isomorphic to $\Z$. Hence, the 
twists are determined by an element of 
$$
H^1(M; {\rm Aut}(\Z)) \cong H^1(M; \Z_2)\;,
$$
since the only nontrivial automorphism of $\Z$ is multiplication by $-1$. 
Twisted integral cohomology may be thought of  as the sheaf cohomology 
$H^\bullet(M; Z)$, taken with respect to the local system $Z$.
%, or defined using a cochain complex. 
The {\v C}ech description goes as follows (see \cite{Fr}). Let $\{ U_\alpha\}$ be an open 
covering of $M$ and $g_{{}_{\alpha \beta}}: U_\alpha \cap U_\beta \to \{\pm1\}$ be 
a cocycle defining the local system $Z$. Then an element of $H^q(M; Z)$ is represented by 
a collection of $q$-cochains $a_i \in Z^q(U_\alpha)$ which satisfy
\(
a_{{}_\beta}=g_{{}_{\alpha \beta}} a_{{}_\alpha} \quad {\rm on} \quad U_{\alpha \beta}
=U_\alpha \cap U_\beta\;.
\label{a}
\)
%We can use any model of cochains, since the group  always acts. 
One can also describe the situation in spaces via Eilenberg-MacLane spaces (see \cite{Fr}). 
More precisely, one uses a model of cochains as maps to Eilenberg-MacLane 
space $K(\Z, q)$. The automorphisms ${\rm Aut}(\Z)\cong\{\pm 1\}$ act on 
$K(\ZZ,q)$ for each $q\geq 0$. For example, $K(\Z, 0)\simeq \Z$ on 
which $-1$ acts by multiplication,
$K(\Z, 1) \simeq S^1$ on which $-1$ acts by reflection. 
The action of ${\rm Aut}(\Z)$ on $K(\Z, q)$ and the cocycle $g_{\alpha \beta}$
gives rise to an associated bundle ${\cal H}^q \to M$ 
with fiber $K(\Z, q)$. 
 Eq. \eqref{a} says that twisted cohomology classes are
represented by sections of  ${\cal H}^q \to M$.
The twisted cohomology group $H^q(M; Z)$ is the set of homotopy classes of 
sections of  ${\cal H}^q \to M$.

\medskip
A third approach uses smooth stacks and unifies the two previous perspectives: namely,
 the {\v C}ech cocycle approach and the approach via Eilenberg-MacLane spaces 
 (see \cite{Cech} \cite{FSS1} \cite{FSS2} for detailed discussions and applications). 
This approach is both extremely general and versatile, and we will rely on it 
when dealing with Deligne cohomology in subsequent sections. 
To that end, let $\cartsp$ denote the category with objects that are convex open subsets of $\RR^n$
with $n\geq 0$ and morphisms that are smooth maps between them. As a coverage on 
this small category we take the covering families to be good open covers  
(i.e. covers with contractible finite intersections). Now let $\BB \ZZ/2$ denote the smooth 
stack on this site obtained by delooping the constant sheaf $\underline{\ZZ/2}$. For example, 
the Dold-Kan image of the sheaf of chain complex $\underline{\ZZ/2}[1]$ can serve as a model.

\medskip 
It turns out that $\BB \ZZ/2$ is a fibrant object in the local projective model structure. 
That is, it is objectwise a Kan-complex and satisfies descent with respect to {\v C}ech hypercovers
 \cite{DI}. Therefore, the connected components of the mapping space $\map(M,\BB \ZZ/2)$ 
 can be calculated by replacing the manifold $M$ by the homotopy colimit over the nerve of 
 an open good open cover $\{U_{\alpha}\}$ (see \cite{Cech} for details). Let $C(\{U_{\alpha}\})$ 
 denote this homotopy colimit. A map
$$
\xymatrix{
M &
\\
\ar[u]^{\simeq}C(\{U_{\alpha}\})\ar[r]^{g} & \BB \ZZ/2
}
$$
determines and is uniquely determined by a {\v C}ech 1-cocycle with coefficients in $\ZZ/2$. 
We, therefore, have an isomorphism
$$
\pi_0\map(C(\{U_{\alpha}\}),\BB \ZZ/2)\cong \check{H}^1(\{U_{\alpha}\};\ZZ/2)\;.
$$

\medskip
The stacky perspective also makes it very transparent how the local system $Z$ and the bundles 
${\cal H}^q$ are related. Indeed, we can similarly define the locally constant, smooth stacks $\BB^q\ZZ$,
 which model higher integral {\v C}ech-cocycles of degree $q$. The action of ${\rm Aut}(\ZZ)\cong \ZZ/2$ 
 on these stacks gives rise to an action groupoid $\BB^q\ZZ/\!/ (\ZZ/2)$, and in turn one can take the homotopy 
 orbit stack associated to this groupoid (see \cite{NSS1} for details). The resulting smooth stack models 
 $\ZZ/2$-bundles with fiber the `stacky Eilenberg-MacLane spaces' $\BB^q\ZZ$. 
Moreover, this bundle is universal in the sense that any such bundle over $M$ fits into a pullback square
\(
\xymatrix{
{\cal H}^q\ar[rr]\ar[d]_-{\pi} && \BB^q\ZZ/\!/(\ZZ/2)=:{\rm Tw}(\BB^q\ZZ)\ar[d]^-{\pi^{\prime}}
\\
M\ar[rr]^-{\eta} && \BB \ZZ/2\;,
}
\label{TwHZ}
\)
defining the twisting stack ${\rm Tw}(\BB^q\ZZ)$. Here the right vertical map is the canonical map 
which projects out $\BB^q\ZZ$ and we have identified the homotopy orbit stack $\ast/\!/(\ZZ/2)$ 
with $\BB \ZZ/2$. Note that we can compute this stack as the homotopy colimit over the nerve of the action
$$
{\rm Tw}(\BB^q\ZZ)\simeq {\rm hocolim}\big\{\xymatrix{
\ar@<-.15cm>[r]\ar@<-.05cm>[r]\ar@<.05cm>[r] \ar@<.15cm>[r]& 
\ZZ/2\times \ZZ/2\times \BB^q\ZZ \ar@<-.1cm>[r]\ar[r]\ar@<.1cm>[r] & 
\ZZ/2\times \BB^q\ZZ\ar@<-.05cm>[r]\ar@<.05cm>[r] & \BB^q\ZZ
}\big\}\;.
$$
\begin{definition}
Define the {\rm $q$-th $\eta$-twisted cohomology} as the connected components of the simplicial 
space of sections $\pi_0\Gamma(M,{\cal H}^q)$. 
\end{definition}

The smooth stacks ${\cal H}^q$ are relatively easy to describe using descent. Indeed, fix a cover 
$\{U_{\alpha}\}$ of $M$ and observe that on each $k$-fold intersection, any map 
$\eta:U_{\alpha_0\hdots \alpha_k}\to \BB\ZZ/2$ factors through the basepoint, since 
$\BB\ZZ/2$ is constant on covers. Consequently, we have 
$$
{\cal H}^q_{U_{\alpha_0\hdots \alpha_k}} \simeq \BB^q\ZZ\times U_{\alpha_0\hdots \alpha_k}\;.
$$ 
Now since $\BB\ZZ/2$ is a 1-type, descent implies that we can compute ${\cal H}^q$ as the coequalizer
$$
{\cal H}^q\simeq {\rm coeq}\Big\{\xymatrix{ \coprod_{\alpha\beta}\BB^q\ZZ\times 
U_{\alpha\beta} \ar@<.05cm>[r]^-{i_{\alpha}\eta_{\alpha\beta}}\ar@<-.05cm>[r]_-{i_{\beta}} &\coprod_{\alpha}\BB^q\ZZ\times U_{\alpha} 
}\Big\}\;,
$$
where $i$ is the inclusion and $\eta_{\alpha\beta}$ is the {\v C}ech cocycle determined by 
$\eta:M\to \BB \ZZ/2$. Then the global sections of ${\cal H}^q$ can be identified with 
a choice of {\v C}ech cocycle $a_{\alpha}$ on each open set $U_{\alpha}$ such that
$
a_{\alpha}=\eta_{\alpha\beta}a_{\beta}
$
on intersections of the cover. We, therefore, recover Freed's description in \cite{Fr}. 

\medskip
We can also recover the local system $Z$ via the sheaf of sections of the bundle ${\cal H}^0$. 
Furthermore, we can see that the sheaf cohomology of $Z$ can be computed as the components of simplicial 
space of sections of the bundles ${\cal H}^q$. We thus have the following equivalent 
characterizations of twisted integral cohomology.

\begin{proposition}
Let $\eta:M\to \BB \ZZ/2$ be a twist for integral cohomology and let $Z$ be the locally constant 
sheaf associated to the twist. Then twisted integral cohomology is given by the isomorphism
$$
H^q_{\eta}(M;\ZZ):=\pi_0\Gamma(M,{\cal H}^q)\cong H^q(M;Z) \;.
$$
\end{proposition}

\section{Deligne cohomology}
%%%%%%%%%%%%%%%%%%

We begin by recalling the definition of smooth Deligne cohomology (see \cite{Del} \cite{Be} \cite{Gi} \cite{Ja} 
\cite{EV} \cite{Ga}). Let ${\cal D}(n)$ denote the sheaf of chain complexes
$$
{\cal D}(n):=\big(
\xymatrix@C=.8cm{ 
\hdots \ar[r] & 0\ar[r] & \underline{\ZZ}\ \ar@{^{(}->}[r] & \Omega^0\ar[r]^{d} & \Omega^1\ar[r]^{d} & \hdots \ar[r]^{d} & \Omega^{n-1}
}\big)\;,
$$
with differential $(n-1)$-forms in degree $0$ and locally constant integer valued functions in degree $n$. For a smooth manifold $M$, the Deligne cohomology group of degree $n$ is defined to be the sheaf cohomology group $\widehat{H}^n(M;\ZZ):=H^0(M;{\cal D}(n))$. These cohomology groups can be explicitly calculated via a {\v C}ech resolution. More precisely, if $\{U_{\alpha}\}$ is a good open cover of $M$, then we can form the {\v C}ech-Deligne double complex
\(\label{double complex}
\xymatrix@=1.8em{
 \underline{\ZZ}(U_{\alpha_0\hdots \alpha_{n}})\ar[r]^-{2\pi i} & \Omega^0(U_{\alpha_0\hdots \alpha_{n}})\ar[r]^-{d} & \Omega^1(U_{\alpha_0\hdots \alpha_{n}})\ar[r]^-{d} & \hdots \ar[r]^-{d} & \Omega^{n-1}(U_{\alpha_0\hdots \alpha_{n}})
\\
\underline{\ZZ}(U_{\alpha_0\hdots \alpha_{n-1}})\ar[r]^-{2\pi i}\ar[u]^-{ (-1)^{n-1}\delta} & \Omega^0(U_{\alpha_0\hdots \alpha_{n-1}})\ar[r]^-{d} \ar[u]^-{(-1)^{n-1}\delta}& \Omega^1(U_{\alpha_0\hdots \alpha_{n-1}})\ar[r]^-{d}\ar[u]^-{(-1)^{n-1}\delta} & \hdots \ar[r]^-{d} & \Omega^{n-1}(U_{\alpha_0\hdots \alpha_{n-1}})\ar[u]^-{(-1)^{n-1}\delta}
\\
\vdots \ar[u]^-{(-1)^{n-2}\delta}& \vdots \ar[u]^-{(-1)^{n-2}\delta}& \vdots \ar[u]^-{(-1)^{n-2}\delta}&  & \vdots \ar[u]^-{(-1)^{n-2}\delta}
\\
 \underline{\ZZ}(U_{\alpha_0\alpha_1})\ar[r]^-{2\pi i}\ar[u]^-{-\delta} & \Omega^0(U_{\alpha_0\alpha_1})\ar[r]^-{d}\ar[u]^-{-\delta} & \Omega^1(U_{\alpha_0\alpha_1})\ar[r]^-{d}\ar[u]^-{-\delta} & \hdots \ar[r]^-{d} & \Omega^{n-1}(U_{\alpha_0\alpha_1})\ar[u]^-{-\delta}
\\
\underline{\ZZ}(U_{\alpha_0})\ar[r]^{2\pi i}\ar[u]^-{\delta} & \Omega^0(U_{\alpha_0})\ar[r]^-{d}\ar[u]^-{\delta} & \Omega^1(U_{\alpha_0})\ar[r]^-{d}\ar[u]^-{\delta} & \hdots \ar[r]^-{d} & \Omega^{n-1}(U_{\alpha_0})\ar[u]^-{\delta}\;,
}
\)
where $U_{\alpha_0\alpha_1\hdots \alpha_k}$ denotes the $k$-fold intersection. The total operator on the double complex is the {\v C}ech-Deligne operator $D:=d+(-1)^p\delta$, where $d$ and $\delta$ 
is the de Rham and {\v C}ech differentials, respectively, acting on elements of degree $p$. 
The sheaf cohomology group $H^0(M;{\cal D}(n))$ can be identified with the group of diagonal elements 
$\alpha_{{}_{k,k}}$ in the double complex which are {\v C}ech-Deligne closed in the sense that 
$(d+(-1)^p\delta)\alpha_{{}_{k,k}}=0$, modulo those which are {\v C}ech-Deligne exact.

\medskip
Deligne cohomology satisfies most of the properties that an ordinary cohomology theory satisfies (such as functoriality and the 
Mayer-Vietoris sequence). However, one needs to be careful when using these properties. For example, 
the form that the Mayer-Vietoris sequence takes is slightly different from what one might expect. The 
following proposition is fairly classical -- a proof for each property can be found in the more modern 
treatment via differential cohomology in \cite{Bu}.
\begin{proposition}[Properties of Deligne cohomology]
Deligne cohomology satisfies the following properties:
\item {\bf (i)} {\rm (Functoriality)} For a smooth map between manifolds $M\to N$, we have 
an induced map
$$
\widehat{H}^n(N;\ZZ)\to \widehat{H}^n(M;\ZZ)\;.
$$
\item {\bf (ii)} {\rm (Additivity)}
For $M=\coprod M_{\alpha}$ a disjoint union of smooth manifolds, we have an isomorphism
$$
\widehat{H}^n(M;\ZZ)\cong\bigoplus_{\alpha}\widehat{H}^n(M_{\alpha};\ZZ)\;.
$$
\item {\bf (iii)} {\rm (Mayer-Vietoris)} For an open cover of $M$ by open smooth manifolds $U$ 
and $V$, we have a sequence
\vspace{-5mm}
$$
\xymatrix{
\hdots \ar[r] & H^{*-2}(U\cap V;\RR/\ZZ) \ar[r] & \widehat{H}^*(M;\ZZ)\ar[r] & 
\widehat{H}^*(U;\ZZ)\oplus \widehat{H}^*(U;\ZZ)  
 \ar@{->} `r/8pt[d] `/12pt[l] `^d[ll]+<-10ex> `^r/8pt[dll] [dll] \\
 & \widehat{H}^*(U\cap V;\ZZ)\ar[r] & H^{*+1}(M;\ZZ)\ar[r] & \hdots \;.
}
$$
\end{proposition}
Note that the Mayer-Vietoris sequence has ordinary integral cohomology on the right and cohomology with 
$\RR/\ZZ$-coefficients on the left. This effect is an artefact of the way the theory is constructed. More 
precisely, the Deligne complex depends on an integer $n$, which indexes the degree of the underlying 
cohomology group. The Mayer-Vietoris sequence comes from the sheaf cohomology of a fixed Deligne 
complex ${\cal D}(n)$ and this is why we see the full differential cohomology group in only one degree. 
In fact, Deligne cohomology is really a mixture of three different cohomology theories (integral, 
$\RR/\ZZ$-coefficients, and de Rham) and captures the interactions between these theories.  
These interactions can be understood via the  ``differential cohomology diagram" 
\(
\label{diamond1}
\xymatrix @C=4pt @!C{
 &\Omega^{*-1}(M)/{\rm im}(d) \ar[rd]^{a}\ar[rr]^{d} & &  \Omega^*_{\rm cl}(M)\ar[rd] &  
\\
H^{*-1}_{\rm dR}(M)\ar[ru]\ar[rd] & & {\widehat{H}^*(M;\ZZ)}\ar[rd]^{I} \ar[ru]^{R}& &  H^*_{\rm dR}(M) \;,
\\
&H^{*-1}(M;\RR/\ZZ) \ar[ru]\ar[rr]^{\beta} & & H^*(M;\ZZ) \ar[ru]^{j} &
}
\) 
where $d$ is the de Rham differential, $R$ is the curvature map, $I$ is the forgetful map, $j$ is 
the rationalization, and $\beta$ is the Beckstein associated with the exponential coefficient sequence. 
This diamond (or hexagon) diagram was originally introduced and emphasized by Simons and 
Sullivan in \cite{SSu} and for more generalized theories, a full characterization via this diamond was 
proved in \cite{BNV}. Parts of it appear in the foundational work of Cheeger and Simons \cite{CS}. 

%\medskip
%Deligne cohomology can be regarded a cochain model for ordinary differential cohomology \cite{Bun}. In the same way that integral cohomology is representable by Eilenberg-Maclane spaces, Deligne cohomology is also representable by certain smooth stacks. More precisely, the moduli stack of higher $n$-bundles with connection $\BB^nU(1)_{\nabla}$ (see \cite{HSS1} for details on this stack) represents Deligne cohomology in degree $n$ and we have an isomorphism
%$$\pi_0\map(M,\BB^nU(1)_{\nabla})\cong \widehat{H}^n(M;\ZZ)\;.$$
%In fact, the computability of Deligne cohomology via {\v C}ech-Deligne cocycles is made even more transparent by again observing that $M$ ought to be replaced by a suitably cofibrant object (for example, the homotopy colimit over its {\v C}ech nerve), in order to compute the homotopy classes of maps on the left hand side. One can show that a map of the form 
%\(
%\xymatrix{
%M &
%\\
%\ar[u]^{\simeq}C(\{U_{\alpha}\})\ar[r]^{g} & \BB^nU(1)_{\nabla}
%}
%\)
%determines and is uniquely determined by a {\v C}ech-Deligne cocycle \cite{FSSt}\cite{Sch}. 

%%%%%%%%%%%%%%%%%%%%%%%
\section{Twists of Deligne cohomology}
\label{Sec twd}
%%%%%%%%%%%%%%%%%%%%%%%

We now discuss the twists for Deligne cohomology. Since Deligne cohomology is a combination of de Rham cohomology  and integral cohomology,  the twists will be some sort of combination of the integral twists 
and twists via differential forms. Note that we have an obvious inclusion of the integral twists
$$
{\rm Aut}(\ZZ)\cong \ZZ/2\into \RR^{\times}\cong {\rm Aut}(\RR)\;.
$$
Delooping this map gives a map between the classifying stacks of units
$$
r:\BB\ZZ/2\into \BB\RR^{\times}\;.
$$
In spaces this map would be an equivalence, but since we are in stacks the geometry prevents this 
map from defining an equivalence of smooth stacks. This map instead simply rationalizes (or realifies) 
the twists and can be viewed as a first approximation to the differential refinement.  

\medskip
The stack $\BB \RR^{\times}$ classifies smooth, locally trivial, real line bundles over a smooth 
manifold and this bundle is determined (up to isomorphism of bundles) by the pullback square
$$
\xymatrix{
{\cal L}_{\eta}\ar[d]\ar[rr] && \RR/\!/\RR^{\times}\ar[d]
\\
M\ar[rr]^-{\eta} && \BB \RR^{\times}\;.
}
$$
Moreover, given a homotopy equivalence between two maps $\eta:M\to \BB \RR^{\times}$ and $\xi:M\to \BB \RR^{\times}$, we have an induced isomorphism of line bundles $\L_{\eta}\to \L_{\xi}$.  This correspondence gives rise to an equivalence of $\infty$-groupoids
$$\map(M,\BB \RR^{\times})\simeq {\bf Line}(M)\;,$$
where on the right we have the $\infty$-groupoid of line bundles on $M$. In a similar way, maps to $\BB \Z/2$ classify $\ZZ$-bundles on $M$ and, given a map $\eta:M\to \BB \ZZ/2$, we have a commutative diagram
$$
\xymatrix@=1.5em{
 & \L_{\eta} \ar[dd]\ar[rr]  && \RR/\!/\RR^{\times}\ar[dd]
\\
Z \ar[ru]^{j}\ar[rr]\ar[dd] && \ZZ/\!/\ZZ/2\ar[ru]\ar[dd] &
 \\
& M\ar[rr]^<<<<<<<<<{i\eta} && \BB \RR^{\times}
\\
M\ar[rr]^-{\eta}\ar@{=}[ru] && \BB \ZZ/2\ar[ru]_{i} &
}
$$
where $Z:={\cal H}^0$ is the $\ZZ$-bundle classified by the map $\eta$. 
The map $j:Z\to \L_{\eta}$ is the bundle map which includes the $\ZZ$-bundle classified by $\eta$ 
as a subbundle of the real line bundle classified by the same twist. 

\medskip
To arrive at the twists of the differential refinement, we need to include the differential form data for the 
rationalization. The crucial ingredient for forming the twisted theory is provided by the \emph{Riemann-Hilbert correspondence}. Fix a flat connection $\nabla$ on a line bundle $\L\to M$ and let $\ker(\nabla)$ denote 
the sheaf of solutions to the parallel transport equation $\nabla s=0$. The flat connection allows us to 
make an identification of sheaves
$$
\ker(\nabla)\simeq \Gamma(-;\L^{\delta})\;,
$$
where the sheaf on the right is the sheaf of sections of the bundle obtained by equipping the fibers of 
$\L\to M$ with the discrete topology and taking transition functions induced by the differences of parallel section on intersections. This identification is at the core of the construction and we 
will revisit it in detail later. We begin by introducing the moduli stack of line bundles with \emph{flat} 
connection. Let $\RR^{\times}:=C^{\infty}(-;\RR^{\times})$ be the sheaf of smooth plots of the punctured real line and let $\Omega^1_{\rm cl}(-)$ be the sheaf of closed 1-forms on the small site of cartesian spaces. The logarithm map $d\log:\RR^{\times}\to \Omega_{\rm cl}^1(-)$ gives 
a smooth action of the sheaf of groups $\RR^{\times}$ on $\Omega_{\rm cl}^1(-)$, 
via the assignment $\omega\mapsto  \omega+d\log(f)$, for sections $\omega\in \Omega_{\rm cl}^1(-)$ 
and $f\in C^{\infty}(-;\RR^{\times})$. 
\begin{definition}
We define {\rm the moduli stack of line bundles with flat  
connection} as the smooth stack 
$$
\flat\BB \RR^{\times}_{\nabla}:=\Omega_{\rm cl}^1(-)\;/\!/\RR^{\times}\;,
$$ 
where the homotopy orbit stack on the right is taken with respect to the action. Alternatively, this smooth
 stack is presented as the image of the positively graded sheaf of chain complexes 
$$
\big(\xymatrix{
\hdots \ar[r] & 0\ar[r] & \RR^{\times}\ar[r]^-{d{\rm log}} & \Omega^1_{\rm cl}(-)
}\big)\;,
$$
under the Dold-Kan functor ${\rm DK}:\chp \to {\rm s}\mathscr{A}{\rm b}\to \sset$. \footnote{Note that we are working over the small site of Cartesian spaces, so stackification of this prestack is not necessary.}
\end{definition}
One can show (see for example \cite{Cech}) that we have an equivalence of $\infty$-groupoids
$$
\map(M,\flat \BB \RR^{\times}_{\nabla})\simeq \flat {\bf Line}(M)\;,
$$
where the infinity groupoid on the right is that of line bundles, equipped with flat connection. 
Let $\underline{\RR}^{\times}$ denote the \emph{constant} smooth stack which associates each element of a good open cover $U_{\alpha}\in \{U_{\alpha}\}$ to the group $\RR^{\times}$.
\footnote{Note that this is different than the smooth stack $\RR^{\times}:=C^{\infty}(-;\RR^{\times})$. 
It is obtained by regarding $\RR^{\times}$ as a discrete group.}
This stack classifies discrete bundles with discrete fiber $\RR^{\delta}$, where here $\RR$ 
is equipped with the discrete topology. By the Poincar\'e Lemma, the morphism
$$
\xymatrix{
\BB \underline{\RR}^{\times}\; \ar@{^{(}->}[r] & \flat \BB \RR^{\times}_{\nabla}\;,
}
$$
from the classifying stack of discrete $\RR$-bundles to 
the classifying stack of $\RR$-bundles with flat connection,  
defines an equivalence on every element of a good open cover. By descent, this implies that we 
have an equivalence of smooth stacks. In summary, we have an induced diagram
\(\label{riemann hilbert}
\xymatrix{
\map(M,\BB \underline{\RR}^{\times})\ar[d] &\ar[l]  \delta{\bf Line}(M)\ar[d]
\\
\map(M,\flat \BB \RR^{\times}_{\nabla})\ar[r] & \flat{\bf Line}(M)
\;,
}
\)
where the right vertical map is defined by the composition and all the maps involved are 
homotopy equivalences of Kan-complexes. It is well-known that the Riemann-Hilbert 
correspondence defines a functor 
$$
{\rm RH}: \flat{\bf Line}(M)\longrightarrow \delta{\bf Line}(M)\;,
$$
which associated each line bundle with flat connection to its corresponding local system. Such a local 
system is equivalently an $\RR^{\delta}$-bundle over $M$ and it is easy to see that ${\rm RH}$ 
defines a homotopy inverse to the right vertical map in \eqref{riemann hilbert}. Restricting to 
Cartesian spaces, this map induces a morphism of smooth stacks
\(\label{RH map}
{\rm RH}:\flat \BB \RR^{\times}_{\nabla} \; \longrightarrow \BB \underline{\RR}^{\times}\;.
\)
This map will help us to identify the twists of the differential refinement, as we shall see. 

\medskip
Essentially the story we are spelling out is that of a twisted de Rham theorem, provided by the 
Riemann-Hilbert correspondence. More concretely, we have the following.
\begin{proposition}\label{twisted de Rham 1}
Let $\nabla$ be a flat connection on a line bundle $\L\to M$ and 
let $\Omega^*(-;\L)$ denote the corresponding de Rham complex with coefficients in 
$\L\to M$ and differential $\nabla$. Let $\L^{\delta}\to M$ be the $\RR^{\delta}$-bundle 
obtained by taking $\L\to M$ to have discrete fibers and constant transition functions 
defined via local parallel sections. Then we have a resolution
$$
j:\L^{\delta}\longrightarrow \Omega^*(-;\L)\;.
$$
In particular, since 
$\Omega^*(-;\L)$ are fine sheaves, this implies that we have an isomorphism
\(\label{twisted de Rham}
H^k(M;\L^{\delta})\cong \frac{{\rm ker}(\nabla:\Omega^k(M;\L)\to 
\Omega^{k+1}(M;\L))}{{\rm im}(\nabla:\Omega^{k-1}(M;\L)\to 
\Omega^{k}(M;\L))}\;.
\)
\end{proposition}
\theproof
Since $\L^{\delta}$ is a locally constant sheaf of vector spaces, the sheaf gluing condition implies 
that for every $U\in \mathcal{O}{\rm pen}(M)$, the sections $\Gamma(U;\L^{\delta})$ appear 
as the kernel
\(\label{exact sequence 1}
\xymatrix{
\Gamma(U;\L^{\delta})\ \ar@{^{(}->}[r] & \prod_{\alpha}\RR
\ar[rr]^-{i_{\alpha\beta}-\eta_{\alpha\beta}i_{\beta\alpha}} && \prod_{\alpha\beta}\RR
}
\)
for some good open cover $\{U_{\alpha}\}$ of $U$, where $\eta_{\alpha\beta}\in \RR^{\times}$ 
are the transition functions of the bundle $\L^{\delta}$. By the fundamental theorem of ODE's, there are nonvanishing local solutions $e_{\alpha}$ to the equation $\nabla(e_{\alpha})=0$ and all such solutions are parametrized by the fiber $\RR$. The $e_{\alpha}$'s define 
local trivializations $\phi_{\alpha}:\RR\times U_{\alpha}\to {\cal L}\vert_{U_{\alpha}}$ by the assignment $(r,x)\mapsto re_{\alpha}(x)$. From these observations, we see that the induced map 
$$j_{\alpha}:\RR\; {\longrightarrow} \; \Gamma(U_{\alpha};\L)\;,$$ 
which maps an element $r\in \RR$ to the corresponding unique solution $e_{\alpha}$, exhibits $\Omega^*(U_{\alpha};\L)$ as a resolution of $\RR$. Now consider the diagram
{\small
\(\label{grid diagram}
\xymatrix@=1.6em{
& 0\ar[d] && 0\ar[d] & 0\ar[d] 
\\
0\ar[r] &\Gamma(U;\L^{\delta})\ \ar[d] \ar[rr]^-{j} && 
\Omega^0(U;\L)\ar[r]^-{\nabla}\ar[d] & \Omega^1(U;\L)\ar[r]\ar[d] & \hdots
\\
0\ar[r] & \prod_{\alpha}\RR \ar[d]^-{i_{\alpha\beta}-\eta_{\alpha\beta}i_{\beta\alpha}} \ar[rr]^-{j_{\alpha} } &&\prod_{\alpha} \Omega^0(U_{\alpha};\L)\ar[d]^-{r_{\alpha\beta}-r_{\beta\alpha}} \ar[r] & \prod_{\alpha}\Omega^1(U_{\alpha};\L)\ar[d]^-{r_{\alpha\beta}-r_{\beta\alpha}} \ar[r] & \hdots 
\\
0\ar[r] &\prod_{\alpha\beta}\RR \ar[rr]^-{j_{\alpha} } \ar[d] &&
 \prod_{\alpha\beta}\Omega^0(U_{\alpha\beta };\L) \ar[r]\ar[d] & 
 \prod_{\alpha\beta}\Omega^1(U_{\alpha\beta };\L) \ar[r]\ar[d] & \hdots
\\
& 0 && 0 & 0
}
\)
where the vertical sequences are short exact and the bottom two horizontal rows are exact. The map $j$ is the map induced by the universal property of the kernel. The diagram commutes as the transition functions $\eta_{\alpha\beta}$ were defined via the local sections $e_{\alpha}$. The top horizontal sequence is natural in $U$ and we have a sequence of sheaves
$$
\xymatrix{
0\ar[r] & \L^{\delta} \ar[r]^-{j} & \Omega^0(-;\L)\ar[r]^-{\nabla}& \Omega^1(-;\L)\ar[r] & \hdots
}
$$
Note that the above diagram \eqref{grid diagram} holds for all refinements of the chosen cover. Using this diagram, along with the fact that the image sheaf is the sheafification of the image in presheaves, \footnote{The sheafification is computed as a limit over refinements of covers and here we know that the sequences are exact.} a quick diagram chase reveals that the above sequence is indeed a resolution of sheaves.
\endofproof

\begin{remark}[Twisted de Rham theorem]
We define the \emph{$\nabla$-twisted de Rham cohomology groups}, $H^*(M;\nabla)$, as the
 quotient on the right hand side of \eqref{twisted de Rham}. Thus, the above proposition states 
 the we have a twisted de Rham isomorphism theorem
$$
H^*(M;\L^{\delta})\cong H^*(M;\nabla)\;.
$$
\end{remark}

Note that we could also define the associated discrete bundle $\L^{\delta}$ as the bundle associated 
to the monodromy representation of the flat connection $\rho:\pi_1(M)\to \RR^{\times}$. For the 
twists of the differential refinement, we need to require that this monodromy representation factors 
through the units of $\ZZ$. That is, we have
$$
\xymatrix{
\rho:\pi_1(M)\ar[r] &  \ZZ/2 \; \ar@{^{(}->}[r] & \RR^{\times}\;.
}
$$
This imposes a restriction on the types of flat connections we can choose on the bundle $\L$. If we start 
with a twist for integral cohomology $\eta:M\to \BB\ZZ/2$ giving the transition functions of a real line 
bundle $\L\to M$, then $\nabla$ must be compatible with this structure. Slight modifications of the 
proof of Proposition \ref{twisted de Rham 1} yield the following.

\begin{proposition}
\label{2nd pro}
Let $\eta:M\to \BB \ZZ/2$ be a twist of integral cohomology and let $\L_{\eta}\to M$ be the real line bundle classified by $\eta$. Let $\rho:\pi_1(M)\to \ZZ/2$ be the map corresponding to the homotopy class of $\eta$ \footnote{The associated map here is provided via the adjunction
$[M,\BB\ZZ/2]\cong [\Pi(M),B\ZZ/2]\cong \hom(\pi_1(M),\ZZ/2)$.} and let $\nabla$ be a flat connection associated to this monodromy representation. Then we have a resolution
$$
j:\L_{\eta}^{\delta}\longrightarrow \Omega^*(-;\L_{\eta})\;,
$$
where $\L_{\eta}^{\delta}$ is the locally constant sheaf obtained via the sheaf of sections of the 
discrete bundle $\L_{\eta}^{\delta}\to M$, obtained by regarding the bundle $\L_{\eta}\to M$ 
as having fiber $\RR$, equipped with the discrete topology. 
\end{proposition} 

Given the information in Proposition \ref{2nd pro}, i.e. a smooth map $\eta:M\to \BB \ZZ/2$ classifying 
a line bundle $\L_{\eta}\to M$, a flat connection $\nabla$ on $\L_{\eta}$ and a resolution 
$j:\L_{\eta}^{\delta}\into \Omega^*(-;\L_{\eta})$, we can define the twisted Deligne complex as follows.
\begin{definition}
\label{def tdc}
{\bf (i)} Given a triple $\nabla:=(\eta,\nabla,j)$ as described above, we define the \emph{twisted Deligne 
complex} as the sheaf of chain complexes on $M$
$$
{\cal D}_{\nabla}(n):=\big(
\xymatrix@C=.8cm{ 
\hdots \ar[r] & 0\ar[r] & Z ~\ar@{^{(}->}[r]^-{j} & \Omega^0(-;\L_{\eta})\ar[r]^-{\nabla} & 
\Omega^1(-;\L_{\eta}) \ar[r]^-{\nabla} 
& \hdots \ar[r]^-{\nabla} & \Omega^{n-1}(-;\L_{\eta})
}\big)\;,
$$
where $\Omega^k(-;\L_{\eta}):=\Omega(-;\Lambda^k(T^*M)\otimes \L_{\eta})$ denotes the 
sheaf of local sections of the bundle and $Z$ is the local system associated to the $\ZZ$-bundle 
classified by $\eta$. 

\item {\bf (ii)} We define the $\nabla$-\emph{twisted Deligne cohomology} 
of $M$ to be the sheaf hypercohomology group
$$
\widehat{H}^n(M;\nabla):=H^0(M;{\cal D}_{\nabla}(n))\;.
$$
\end{definition}

The twists for the Deligne complex can be organized into a smooth stack themselves. 

\begin{definition}\label{deligne twists}
We define the {\rm stack of twists for the Deligne complex} as the $(\infty,1)$-pullback
$$
\xymatrix{
\BB (\ZZ/2)_{\nabla}\ar[rr]\ar[d] && \flat (\BB \RR^{\times}_{\nabla})\ar[d]^-{\rm RH}
\\
\BB \ZZ/2\ar[rr]^-{r} &&  \BB \underline{\RR}^{\times}
\;.
}
$$
\end{definition}

This smooth stack indeed defines the necessary information.
\begin{proposition}
A map $M\to \BB(\ZZ/2)_{\nabla}$ determines and is uniquely determined by the following data:
\begin{enumerate}
\item A discrete $\RR^{\delta}$-bundle $\L^{\delta}_{\eta}\to M$ classified by a map 
$\eta:M\to \BB \ZZ/2$ and a $\ZZ$-subbundle $Z\to M$.
\item A flat connection $\nabla$ on a line bundle $\L\to M$.
\item An isomorphism of local systems 
$j:\L^{\delta}_{\eta}\overset{\cong}{\longrightarrow} {\rm ker}(\nabla)$ giving rise to 
a resolution $j:\L^{\delta}_{\eta}\to \Omega^*(-;\L)$.
\end{enumerate}
\end{proposition}
\theproof
A map $M\to \BB(\ZZ/2)_{\nabla}$ can be identified with a pair of maps $\eta:M\to \BB \ZZ/2$ 
and $\nabla:M\to \flat(\BB \underline{\RR}^{\times}_{\nabla})$ such that ${\rm RH}(\nabla)$ 
and $r(\eta)$ are connected by an edge in 
$\map(M,\BB \RR^{\times})$. 
An edge in this mapping space can be identified with an isomorphism of corresponding 
$\RR^{\delta}$-bundles. Thus, we must have an isomorphism of corresponding local systems
$$
j:\L_{\eta}^{\delta}\longrightarrow {\rm ker}(\nabla)\;.
$$
By definition, $\Omega^*(-;\L)$ resolves ${\rm ker}(\nabla)$ and therefore the isomorphism 
gives rise to the desired resolution.\endofproof
\begin{remark}Since the map $RH$ in definition \ref{deligne twists} an equivalence, the induced map $\BB(\ZZ/2)_{\nabla}\to \BB \ZZ/2$ is also an equivalence.  Thus, up canonical equivalence, there is a unique differential refinement of any topological twist $\eta$.
\end{remark}

We illustrate the definition with the following example. 

\begin{example}[Punctured complex plane]
\label{Ex1}
Let $M$ be the punctured complex plane $\CC{-}\{0\}$. There are two isomorphism classes of real line 
bundles on $\CC{-}\{0\}$, classified by $H^1(\CC{-}\{0\};\ZZ/2)\cong \ZZ/2$: the trivial bundle and the 
M\"obius bundle. Let $\mathcal{L}\to \CC{-}\{0\}$ denote the M\"obius bundle and $\nabla$ be a flat connection compatible with the monodromy representation defined by sending $1\in \ZZ\cong \pi_1(S^1)$ to $-1\in \ZZ/2$. Notice that this representation also defines the principal $\ZZ/2$-bundle associated to the M\"obius bundle over $\CC{-}\{0\}$. Consider the open cover $\{U,V\}$ of $\CC{-}\{0\}$ obtained by removing the rays $x>0$ and $x<0$, where $z=x+iy$. Let $s_{U}$ be the local section traversing one edge of M\"obius strip on $U$ and $s_{V}$ be the local section traversing the same edge on $V$. Then $s_{U}$ and $s_{V}$ define local trivializations, in which $\nabla=d$.
In this case, the twisted Deligne complex takes the form 
 $$
{\cal D}_{\nabla}(n):=\big(
\xymatrix@C=.8cm{ 
\hdots \ar[r] & 0\ar[r] & Z ~\ar@{^{(}->}[r] & 
\Omega^0(-;\L)\ar[r]^-{\nabla} & \Omega^1(-;\L) \ar[r]^-{\nabla} 
& \hdots \ar[r]^-{\nabla} & \Omega^{n-1}(-;\L)
}\big)\;,
$$
where $Z$ is the sheaf of sections of the $\ZZ$-subbundle of the M\"obius bundle. Of course, for dimension reasons, we only need to consider ${\cal D}_{\omega}(n)$ up to degree $n=2$. Locally this complex is \emph{isomorphic} to the untwisted Deligne complex and the isomorphism is defined by $s_{U}$ and $s_{V}$. For a general nonvanishing local sections $\sigma_{U}$ and $\sigma_{V}$, the complex will be isomorphic (over $U$ for example) to  
 $$
{\cal D}_{\nabla}(n):=\big(
\xymatrix@C=.8cm{ 
\hdots \ar[r] & 0\ar[r] & \underline{\ZZ} ~\ar@{^{(}->}[r] & 
\Omega^0(-)\ar[r]^-{d+df_{U}} & \Omega^1(-) \ar[r]^-{d+df_{U}} 
& \hdots \ar[r]^-{d+df_{U}} & \Omega^{n-1}(-)
}\big)\;,
$$
with $f_{U}$ a smooth function such that $f_{U}s_{U}=\sigma_{U}$. 
%
%Indeed, the monodromy representation $\rho:\pi_1(\CC-\{0\})\cong \ZZ\to \ZZ/2$, defined by the connection gives the correct associated principal bundle. The holonomy subbundle in the degenerate case, i.e. $k=0$,
% gives the local system $\underline{\ZZ}$. In the nondegenerate case, we get nested helices, with associated 
% local coefficient systems $Z_{\omega_k}$. Since in each case the bundle is trivial, we arrive at the twisted 
% Deligne complex
% $$
%{\cal D}_{\omega_k}(n):=\big(
%\xymatrix@C=.8cm{ 
%\hdots \ar[r] & 0\ar[r] & Z_{\omega_k} ~\ar@{^{(}->}[r] & \Omega^0\ar[r]^-{d+\omega_k} & \Omega^1 \ar[r]^-{d+\omega_k} 
%& \hdots \ar[r]^-{d+\omega_k} & \Omega^{n-1}
%}\big)\;.
%$$
%Then, when $k=0$ we arrive at the usual Deligne complex. Note that even though each $\omega_k$ is equivalent to the zero connection, the corresponding twisted Deligne complexes are not quasi-isomorphic! Thus, the set of twists are actually more refined then equivalence classes of flat connections.
\end{example}

\begin{remark}
Note that there is a canonical map from the twisted Deligne cohomology groups to the twisted de Rham cohomology groups. Indeed, for any line bundle ${\cal L}\to M$ and connection $\nabla$, we have a morphism of complexes 
 \(\label{map of chain}
\xymatrix@C=.8cm{ 
\hdots \ar[r] & 0\ar[r] & Z ~\ar@{^{(}->}[r]^-{j}\ar[d] & \Omega^0(-;{\cal L}) \ar[r]^-{\nabla}\ar[d] & \Omega^1(-;{\cal L}) \ar[r]^-{\nabla}\ar[d]
& \hdots \ar[r]^-{\nabla} & \Omega^{n-1}(-;{\cal L})\ar[d]^{\nabla}
\\
\hdots \ar[r] & 0\ar[r] & 0 \ar[r] & 0\ar[r] & 0 \ar[r] 
& \hdots \ar[r] & \Omega_{\rm fl}^{n}(-;{\cal L})
\;,
}
\)
since $\nabla^2=0$. Here, $\Omega_{\rm fl}^{n}(-;{\cal L})$ denotes the subsheaf 
${\rm ker}\big(\nabla:\Omega^{n}(-;{\cal L})\to \Omega^{n+1}(-;{\cal L})\big)$.
\end{remark}

\begin{remark} 
The map \eqref{map of chain} then induces a map
$R:\widehat{H}^n(M;\nabla)\to H^n_{\rm dR}(M;\nabla)$.
Note that $R$ is in fact natural in $M$ in the following sense. Let $f:M\to N$ be a smooth map and 
fix a twist $\nabla$ on $N$. Let $f^*(\nabla)$ be the pullback of the connection on the line bundle 
$f^*(\L)$. Via functoriality of sheaf cohomology, we get an induced map
$$
f^*:\widehat{H}^n(N;\nabla)\longrightarrow \widehat{H}^n(M;f^*(\nabla))\;.
$$ 
Then we have a commutative diagram
$$
\xymatrix{
\widehat{H}^n(N;\nabla)\ar[r]^-{R}\ar[d]_-{f^*} & H^n_{\rm dR}(N;\nabla)\ar[d]^-{f^*}
 \\
 \widehat{H}^n(M;f^*(\nabla))\ar[r]^-{R} & H^n_{\rm dR}(M;f^*(\nabla))
 \;.
 }
 $$
 \end{remark}
%The following simple lemma will be useful in characterizing the twists 
%\begin{lemma}
%Let $\eta:M\to \BB\ZZ/2$ be an integral twist with associated bundle ${\L}_{\eta}$. Every connection $\nabla$, compatible with the $\ZZ/2$-structure, is of the form $d+\omega$, for some globally defined 1-form $\omega$. 
%\end{lemma}
%\theproof
%Write $\nabla=d+\omega_{\alpha}$ on open sets $U_{\alpha}$ covering $M$. The map $\eta:M\to \BB\ZZ/2$ specifies transition functions $\eta_{\alpha\beta}:U_{\alpha\beta}\to \ZZ/2$. Since these functions are $\ZZ/2$-valued, $\nabla$ transforms as
%$$\omega_{\alpha}\mapsto \omega_{\alpha}+d\eta_{\alpha\beta}=\omega_{\alpha}$$
%on intersections. By the sheaf condition, $\omega_{\alpha}$ are the restriciton of a globally defined form $\omega$. 
%\endofproof

The next proposition shows that twisted Deligne cohomology indeed reduces to ordinary cohomology when 
the twist is trivial.
When $\nabla\simeq 0$ the local systems $Z$ and $R$ trivialize: $Z\leadsto \underline{\ZZ}$ and $R\leadsto \underline{\RR}$. Moreover, $\omega=0$, and we recover the usual calculation for cohomology with 
coefficients in the Deligne complex.

\begin{proposition}\label{trivialization of twists}
Let $\nabla:M\to \BB(\ZZ/2)_{\nabla}$ be a twist of Deligne cohomology which is trivial in the sense that 
$\nabla$ factors through the basepoint 
$0:\ast\to \BB (\ZZ/2)_{\nabla}$ up to homotopy 
equivalence. Then we have a natural isomorphism of functors
$$
\widehat{H}^n(-;\ZZ)\cong \widehat{H}^n(-;\nabla)\;.
$$
\end{proposition}
\theproof
If $\nabla\simeq 0$, then in particular $\eta\simeq 0$ and the bundle $\L\to M$ is trivializable as a bundle with flat connection. Moreover, the homotopy gives rise to a preferred choice of trivialization for each structure. These trivializations gives rise to a quasi-isomorphism 
of complexes
 $$
\xymatrix@C=.8cm{ 
\hdots \ar[r] & 0\ar[r] & Z ~\ar@{^{(}->}[r]^-{j}\ar[d] & \Omega^0(-;{\cal L}) \ar[r]^-{\nabla}\ar[d] & \Omega^1(-;{\cal L}) \ar[r]^-{\nabla}\ar[d]
& \hdots \ar[r]^-{d} & \Omega^{n-1}(-;{\cal L})\ar[d]
\\
\hdots \ar[r] & 0\ar[r] & \underline{\ZZ} ~\ar@{^{(}->}[r]^-{i} & \Omega^0(-) \ar[r]^-{d} & \Omega^1(-) \ar[r]^-{d}
& \hdots \ar[r]^-{d} & \Omega^{n-1}(-)
\;,
}
$$
and the claim follows.
\endofproof

%%%%%%%%%%%%%%%%%%%%%%%%%%
\section{Properties of twisted Deligne cohomology}
\label{Sec Prop}
%%%%%%%%%%%%%%%%%%%%%%%%%%

In this section we discuss the properties of basic twisted Deligne cohomology. Several of these properties 
have familiar counterparts in ordinary cohomology, while others have properties which are analogous to 
ordinary differential cohomology. We start with calculating the sheaf cohomology groups with coefficients 
the local systems associated with the twists.

\begin{lemma}\label{cohomology calc}
Let $\nabla:M\to \BB(\ZZ/2)_{\nabla}$ be a twist for the Deligne complex on a smooth manifold $M$ and let $\eta:M\to \BB \ZZ/2$ denote the underlying topological twist. The sheaf cohomology groups of degree $k\neq 0$ of 
${\cal D}_{\nabla}(n)$ (see Def. \ref{def tdc}) are given by
$$
H^k(M;{\cal D}_{\nabla}(n))\cong\left\{\begin{array}{lcc}
H^{n+k}(M;Z), && k>0
\\
\\
H^{n+k-1}(M;\L^{\delta}/Z), && k<0.
\end{array}\right.
$$
\end{lemma}
\theproof 
In what follows, the complex $\Omega(-;\L)^*$ is equipped with the differential $\nabla$. The complex ${\cal D}_{\nabla}(n)$ 
is quasi-isomorphic to the shifted cone ${\rm cone}\big(Z\oplus \tau_{\leq 0}\Omega(-;\L)*[n]\to \Omega(-;\L)^*\big)[-1]$, where the 
morphism is given by the assignment $(z,x)\mapsto j(z)-x$. Here $\tau_{\leq 0}$ is the truncation functor to positive degrees. 
Thus, we have an exact triangle
\(
\hspace{-.9cm}
\label{triangle}
\xymatrix{
& {\rm cone}\big(Z\oplus \tau_{\geq 0}\Omega(-;\L)^*[n]\to \Omega(-;\L)^*[n]\big)[-1]
\ar[r] & \tau_{\geq 0}\Omega(-;\L)^*[n]\oplus Z 
 \ar@{->} `r/8pt[d] `/12pt[l] `^d[l]+<-10ex> `^r/8pt[dl] [dl] \\
& ~~\Omega^*(-;\L)[n]  \ar[r]& {\rm cone}\big(Z\oplus \tau_{\geq 0}\Omega(-;\L)*[n]\to \Omega(-;\L)^*[n]\big)\;.
 }
\)
Note that we also have a short exact sequence
\(\label{short sequence 1}
0
\longrightarrow
 {\rm cone}\big(Z\to \Omega(-;\L)^*[n]\big)[-1]
\longrightarrow
{\cal D}_{\nabla}(n)
 \longrightarrow
 {\rm cone}\big(\tau_{\geq 0}\Omega(-;\L)^*[n]\to 0\big)[-1]
 \longrightarrow
  0
\)
induced by the map $j$ and the projection. Since $\Omega^*(-;\L)$ resolves the locally constant sheaf 
$\L^{\delta}$, the mapping cone on the left is quasi-isomorphic to $\L^{\delta}/Z[n-1]$. The mapping 
cone on the right has trivial sheaf cohomology in negative degrees. The long exact sequence in sheaf
 cohomology gives an isomorphism
%\vspace{-3mm}
$$
H^k(M;{\cal D}_{\nabla}(n)) \cong H^k(M;\L^{\delta}/Z[n-1])
 \cong H^{n+k-1}(M;\L^{\delta}/Z)\;,
$$
for $k<0$. For $k>0$, the long exact sequence obtained from the exact triangle
\eqref{triangle}
 gives the desired isomorphism.
\endofproof

Abstractly, a twisted differential cohomology theory should satisfy certain axioms
and properties.  We now  verify an explicit  diagrammatic characterization of twisted Deligne 
cohomology, refining diagram \eqref{diamond1}.

\begin{proposition}[Twisted Deligne cohomology diamond]
\label{tw diamond}
Let $M$ be a smooth manifold and let $\nabla:M\to \BB (\ZZ/2)_{\nabla}$ be a twist 
for the Deligne complex on $M$. 
Then the twisted Deligne cohomology groups fit into the diamond diagram
$$
\xymatrix @C=4pt @!C{
 &\Omega(M;\L)^{*-1}/{\rm im}(\nabla) \ar[rd]^{a}\ar[rr]^{\nabla} & &  
 \Omega(M;\L)^*_{\rm fl}\ar[rd] &  
\\
H_{\rm dR}^{*-1}(M;\nabla)\ar[ru]\ar[rd] & & 
{\widehat{H}^*(M;\nabla)}\ar[rd]^{I} \ar[ru]^{R}& &  H^*_{\rm dR}(M;\nabla) \;.
\\
&H^{*-1}(M;\L^{\delta}/Z) \ar[ru]\ar[rr]^{\beta} & & H^*(M;Z) \ar[ru]^{j} &
}
$$
\end{proposition}
\theproof
The right square follows from the sequence \eqref{triangle}, after passing to sheaf cohomology. 
The two diagonal sequences are induced from the two short exact sequences, 
$$
0
\longrightarrow
 \tau_{<0}\Omega(-;\L)^*[n]
 \longrightarrow
{\cal D}_{\nabla}(n)
  \longrightarrow
   Z[n]
   \longrightarrow
    0
$$
and \eqref{short sequence 1}. The commutativity of the top part of the diagram follows from the 
web of short exact sequences
{
\footnotesize 
$$
\hspace{-.7cm}
\xymatrix@C=.7em{
{\rm cone}\big(0\to \tau_{\geq 0}\Omega(-;\L)^*[n]\big)[-1] \ar[r]\ar[d] &
 {\rm cone}\big(Z\to \Omega(-;\L)^*[n]\big)[-1] \ar[r]\ar[d] &
  {\rm cone}\big(Z\to \tau_{<0}\Omega(-;\L)^*[n]\big)[-1] \ar[d]
\\
{\rm cone}\big(\tau_{\geq 0}\Omega(-;\L)^*[n]\to \tau_{\geq 0}\Omega(-;\L)^*[n]\big)[-1] \ar[r]\ar[d] &
{\rm cone}\big(Z\oplus \tau_{\leq 0}\Omega(-;\L)^*[n]\to \Omega(-;\L)^*\big)[-1]\ar[r]^-{\pi}_-{\simeq} \ar[d]^{R}  
& {\cal D}_{\nabla}(n) \ar[d]
\\
{\rm cone}\big(\tau_{\geq 0}\Omega(-;\L)^*[n]\to 0\big)[-1]\ar@{=}[r] & 
{\rm cone}\big(\tau_{\geq 0}\Omega(-;\L)^*[n]\to 0\big)[-1] \ar[r] & 0 
}
$$}
along with the fact that the connecting homomorphisms $\delta_v$ and $\delta_{h}$ for the 
vertical and horizontal sequences, respectively, obey \footnote{It is a straightforward exercise 
to show that this follows in general for any such web of short exact sequences.} 
$$
\delta_vR=-\delta_{h}\pi\;.
$$
For commutativity of the bottom part of the diagram, consider the long exact sequence associated to 
the cone ${\rm cone}(Z\to \Omega^*(-;\L)[n])[-1]$. This sequence is just the Bockstein sequence associated 
to the short exact sequence
$$
0\longrightarrow Z \longrightarrow \L^{\delta} \longrightarrow 
\L^{\delta}/Z\longrightarrow 0\;,
$$
shifted down by 1. Taking long exact sequences associated to cones, the first map in 
sequence \eqref{short sequence 1} induces a commutative diagram
$$
\hspace{-.2cm}
\xymatrix@C=1.7em{
H^{n-1}(M;Z)\ar[r]\ar@{=}[d] & H^{n-1}(M;\L^{\delta})\ar[r]\ar[d]^{\cong} & 
H^{n-1}(M;\L^{\delta}/Z)\ar[r]^{-\beta}\ar[d] & H^n(M;Z) \ar[r]\ar[d] & H^n(M;\L^{\delta})\ar[d]^{\cong}
\\
H^{n-1}(M;Z)\ar[r] & H^{n-1}_{\rm dR}(M;\nabla)\ar[r] & \widehat{H}^n(M;\nabla) \ar[r]^-{{I}\oplus {R}} & 
H^n(M;Z)\oplus \Gamma^n_{\rm fl}(M;\L) \ar[r] & H_{\rm dR}^n(M;\nabla)\;. 
}
$$
Therefore, the result follows. 
\endofproof

Besides the differential cohomology diamond, we also have a Mayer-Vietoris sequence at our disposal. Combined with  
Lemma \ref{cohomology calc}, this takes the following form.

\begin{proposition} [Mayer-Vietoris for twisted Deligne cohomology]
Let $M$ be a smooth manifold with open cover $\{U,V\}$ and let $\nabla:M\to \BB (\ZZ/2)_{\nabla}$ be a twist with underlying topological twist $\eta$. There is a Mayer-Vietoris type sequence
\vspace{-1mm}
$$
\xymatrix{
\hdots \ar[r] & H^{*-2}(U\cap V;\L^{\delta}/Z) \ar[r] & \widehat{H}^*(M;\nabla)\ar[r] & 
\widehat{H}^*(U;\nabla)\oplus \widehat{H}^*(V;\nabla)  
 \ar@{-<} `r/8pt[d] `/12pt[l] `^d[ll]+<-10ex> `^r/8pt[dll] [dll] \\
 & ~\widehat{H}^*(U\cap V;\nabla)\ar[r] & H^{*+1}(M;Z)\ar[r] & \hdots \;.
}
$$
\end{proposition}
\theproof
The Mayer-Vietoris sequence holds for sheaf cohomology, hence for ${\cal D}_{\nabla}(n)$. The claim then follows from the characterization of Lemma 
\ref{cohomology calc}.
\endofproof

%%%%%%%%%%%%%%%%%%%
\section{Computations and examples}
\label{Sec Ex}
%%%%%%%%%%%%%%%%%%%

In this section we compute the twisted Deligne cohomology for various spaces, illustrating the 
constructions and computational techniques developed earlier. We start with the 
simplest case. 

\begin{example}[Twisted Deligne cohomology of $\RR^n$]
\label{Ex Rn}
Since $H^1(\RR^n;\ZZ/2)\cong 0$, every real line bundle is trivializable over $\RR^n$. A flat 
connection on a trivial line bundle is simply a closed differential 1-form $\omega$. Since the de 
Rham cohomology of $\RR^n$ is also trivial, $\omega=d\beta$, and multiplication by the 
exponential map gives a quasi-isomorphism of complexes
$$
e^{\beta}\times:(\Omega^*(-),d+d\beta)\longrightarrow (\Omega^*(-),d)\;.
$$  
Thus, we see that the twisted Deligne cohomology groups reduce to the ordinary Deligne 
cohomology groups. These in turn are easily computed as
$$\widehat{H}^0(\RR^n;\ZZ)\cong \ZZ\;,
\qquad 
\widehat{H}^1(\RR^n;\ZZ)\cong C^{\infty}(\RR^n,\RR/\ZZ)\;,
\qquad
\text{and}
\qquad
\widehat{H}^k(\RR^n;\ZZ)\cong \Omega^{k-1}(\RR^n)/{\rm im}(d)\;.
$$
\end{example}

Similar effects occur for the punctured Euclidean space. 

\begin{example}[Twisted Deligne cohomology of 
$\RR^n{-}\{0\}$, $n>1$]
\label{Ex pun}
For $n>1$, the punctured real $n$-space $\RR^n{-}\{0\}$ is simply connected, we have $H^1(\RR^n{-}\{0\};\ZZ/2)\cong 0$. Therefore, every real line bundle over 
$\RR^n{-}\{0\}$ is trivializable. Similarly, every closed differential
 1-form is exact. As in the calculation of the differential cohomology of $\RR^n$ (Example \ref{Ex Rn}), 
 it follows that the twisted Deligne cohomology groups reduce to the ordinary Deligne cohomology
  groups. These are readily calculated via the diamond (Prop. \ref{tw diamond}) as
$$
\widehat{H}^0(\RR^n{-}\{0\};\ZZ)\cong \ZZ\
\qquad 
\text{and} 
\qquad
\widehat{H}^1(\RR^n{-}\{0\};\ZZ)\cong C^{\infty}(\RR^n{-}\{0\},\RR/\ZZ)\;.
$$
For $k\neq n-1,n$, we have
$$
\widehat{H}^k(\RR^n{-}\{0\};\ZZ)\cong \Omega^{k-1}(\RR^n{-}\{0\})/{\rm im}(d)\;.
$$
For $k=n-1$, we have the identification, via the Hodge decomposition, 
$$
\widehat{H}^{n-1}(\RR^n{-}\{0\};\ZZ)\cong \Omega^{n-2}(\RR^n{-}\{0\})/{\rm im}(d)
\oplus \langle \omega \rangle\;,
$$
where $\langle \omega \rangle$ is the $\ZZ$-linear span of a normalized harmonic $(n-1)$-form, restricting to a volume form $(n-1)$-sphere. Note also that the identification depends on a choice of metric. For $k=n$ and a choice of metric, the Hodge decomposition now gives
$$
\widehat{H}^{n}(\RR^n{-}\{0\};\ZZ)\cong d^{\dagger}\Omega^{n}(\RR^{n}{-}\{0\})\oplus \RR/\ZZ\;,
$$
where the copy of $\RR/\ZZ$ is identified with the group $\langle \omega\rangle_{\RR} /\langle \omega\rangle_{\ZZ}$, with $\omega$ the harmonic form extending the normalized volume form of the $(n-1)$-sphere and the subscripts indicate that we are taking the $\RR$ and $\ZZ$-linear spans, respectively.
\end{example}

Passing to the complex setting allows for some additional information in the twists. 

\begin{example}[Twisted Deligne cohomology of the punctured complex plane]
\label{Ex2}
Let $\nabla$ be the flat connection on the M\"obius bundle $\L\to \CC{-}\{0\}$ as in example \ref{Ex pun}. Recall that we have nontrivial monodromy arising from the map 
$$
\xymatrix{
\pi_1(\CC{-}\{0\})\cong \ZZ\ar[r] & \ZZ/2 \; \ar@{^{(}->}[r] & \RR^{\times}
}
\;,
$$
sending $1\mapsto -1$. Let $Z$ denote the local system corresponding to the $\ZZ$-subbundle of the M\"obius bundle. 
We can compute the cohomology with coefficients in $Z$ via {\v C}ech cohomology, as follows.
 Consider the covering of $\CC{-}\{0\}$ by three open sets as in the following picture:
 
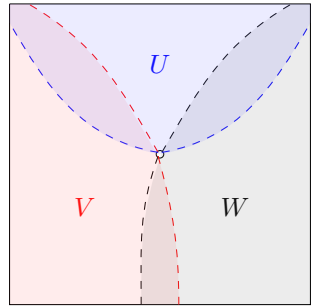
\begin{center}
\begin{tikzpicture}
\draw (-2,2) -- (2,2) -- (2,-2) -- (-2,-2) -- (-2,2);
\draw (0,0) circle (.05cm);
%open sets
\node at (0,1.2)  {${\color{blue} U}$};
\node at (-1,-.7)  {${\color{red} V}$};
\node at (1,-.7)  {$W$};

\draw[dashed,red] (-.025,.025) to[out=120,in=-30] (-1.75,2);
\draw[dashed,blue] (-.025,.025) to[out=170,in=-60] (-2,1.75);

\draw[dashed,black] (.025,.025) to[out=60,in=-150] (1.75,2);
\draw[dashed,blue] (.025,.025) to[out=10,in=-120](2,1.75);

\draw[dashed,black] (-.025,-.025) to[out=-110,in=90] (-.25,-2);
\draw[dashed,red] (-.025,-.025) to[out=-75,in=90] (.25,-2);

\path[fill=red!50,opacity=.15] (-1.75,2)  to[out=-30,in=120] (-.025,-.025) to[out=-75,in=90] (.25,-2) to (-2,-2) to (-2,2) to (-1.75,2);

\path[fill=black!50,opacity=.15]  (1.75,2) to[out=-150,in=60] (.025,.025) to[out=-110,in=90] (-.25,-2) to (2,-2) to (2,2) to (1.75,2);

\path[fill=blue!50,opacity=.15]   (-2,1.75)  to[out=-60,in=170] (0,0) to[out=10,in=-120](2,1.75) to (2,2) to (-2,2) to (-2,1.75);

\end{tikzpicture}
\end{center}
where the boundary of each of the open sets $V,U$ and $W$ is colored red, blue and black, respectively. Note that, because the plane is punctured, there are no three-fold intersections of the open sets. The two-fold intersections are the wedge 
regions in between the colors. Now we can choose a cocycle representative for the twist defined by the assignment 
of integers $(-1,-1,-1)$ on the double intersections depicted above. In integral cohomology, the cocycle $(-1,-1,-1)$ is a representative for the 
generator of $H^1(\CC{-}\{0\};\ZZ)\cong \ZZ$ and is \emph{not} an integral 
{\v C}ech coboundary. However, when viewed as a {\v C}ech cocycle in the sheaf $Z$, $2(-1,-1,-1)=(-2,-2,-2)$ is a coboundary. Indeed, the transition functions for the bundle modify the restriction maps for $Z$ and we are reduced to showing that the system of equations
$$
n_{U}-(-n_{V})=-2,\qquad n_{U}-(-n_{W})=-2,\qquad n_{V}-(-n_{W})=-2
$$
has a solution. But this is easy to see, for example $n_{U}=-1$, $n_{V}=-1$, $n_{W}=-1$. The since the analogous system for the cocycle $(-1,-1,-1)$ has no solutions, we conclude that
$$
H^1(\CC{-}\{0\};Z)\cong \ZZ/2\;.
$$
To see what a global section of $Z$ looks like, we attempt to find solutions the analogous equations
$$
n_{U}-(-n_{V})=0,\qquad  n_{U}-(-n_{W})=0,\qquad n_{V}-(-n_{W})=0\;.
$$
But these imply that $n_{V}=-n_{W}$, $n_{U}=n_{W}$ and $2n_{U}=0$. Hence, $n_{U}=n_{V}=n_{W}=0$ and so
$$
H^0(\CC{-}\{0\};Z)\cong \Gamma(\CC{-}\{0\},Z)\cong 0\;.
$$
The twisted de Rham cohomology is now easy to compute from the twisted de Rham theorem. Indeed, the calculations in {\v C}ech cohomology apply equally well to the discrete bundle $\L^{\delta}$. The presence of 2-torsion in degree 1 kills $H^1_{\rm dR}(\CC{-}\{0\};\nabla)$ and we have the identifications
$$
H_{\rm dR}^0(\CC{-}\{0\};\nabla)\cong H_{\rm dR}^1(\CC{-}\{0\};\nabla)\cong 0\;.
$$
From the differential cohomology diamond diagram (Prop. \ref{tw diamond}), we have $\widehat{H}^0(\CC{-}\{0\};\nabla)\cong \widehat{H}^0(\CC{-}\{0\};Z)\cong \ZZ$. Since the global sections of $\L$ are divisible as an abelian group, they form an injective module. By the differential cohomology diamond diagram (Prop. \ref{tw diamond}), we conclude that 
$$
\widehat{H}^1(\CC{-}\{0\};\nabla)\cong \ZZ/2\oplus \Gamma(\CC{-}\{0\};\L)
\qquad \text{and} \qquad 
\widehat{H}^0(\CC{-}\{0\};\nabla)\cong 0\;.
$$

\medskip
\end{example}
\begin{example}[Orientation line bundle]
Let $M$ be a closed, smooth manifold of dimension $n$, with orientation bundle $\Lambda^nM\to M$. 
If $M$ is simply-connected then $\Lambda^nM\to M$ is trivializable and every closed 1-form is exact. 
In this case, the twisted Deligne cohomology groups reduce to the usual Deligne cohomology groups. 

Let $M$ be a non-orientable manifold and equip $M$ with a Riemannian metric so that $TM\to M$ has 
orthogonal structure. Consider the Levi-Civita connection $\nabla$ on $TM$. Taking the determinant of the transition functions of $TM$ and the trace of the connection gives the orientation bundle
 $\Lambda^nM\to M$, equipped with the zero connection. 

Let $U_{\alpha}$ be a local chart of $M$, with coordinates $\{x_i\}_{i=1}^n$. Write a differential $n$-form 
on a patch $U_{\alpha}$ as $f_{\alpha}dx_1\wedge dx_2\wedge \hdots \wedge dx_n$, with $f\in C^{\infty}(U_{\alpha};\RR)$. From the definitions, we see that a flat section of this bundle is locally of the above form with $f_{\alpha}\equiv C\in \RR$. Under orthonormal coordinate transformations, we see that these constants differ by the determinant ${\rm det}(g_{\alpha\beta})=\pm 1$. Thus, the sheaf of local sections can be regarded as sections of the bundle $\Lambda^nM^{\delta}\to M$, where the $\delta$ indicates that we have taken the fibers $\RR$ to have the discrete topology. The local system $Z\subset \Lambda^nM^{\delta}$ in this case is the sheaf of sections which are locally of the form $n_{\alpha}dx_1\wedge dx_2\wedge \hdots \wedge dx_n$, with 
$n_{\alpha}\in \ZZ$.

In this case, the $\nabla=d_{\Lambda^nM}$-twisted Deligne complex is given by the hypercohomology of the complex
$$\xymatrix{
Z\ \ar@{^{(}->}[r] & \Omega^0(-;\Lambda^nM)\ar[r]^-{d} & \Omega^1(-;\Lambda^nM)\ar[r]^-{d} & \hdots \ar[r]^-{d} & \Omega^{n-1}(-;\Lambda^nM)
}\;.
$$
\end{example}

\medskip
\begin{example}[Twisted Deligne cohomology of real projective space]
\label{ex3}
 Let $M=\RR P^{n}$ be the $n$-dimensional real projective space. The restriction of the first Stiefel-Whitney 
 class $w_1\in H^1(\RR P^{\infty};\ZZ/2)\simeq \ZZ/2$ to $\RR P^{n}$ classifies the tautological bundle on 
 $\RR P^n$. This tautological bundle is a special case of the previous example. Hence, the $w_1$-twisted 
 Deligne complex reduces to
$$
{\cal D}_{w_1}(k)=\Big(
\xymatrix{
Z\ \ar@{^{(}->}[r] &\Gamma^0(-; \L_{w_1}) \ar[r]^-{d} & \Gamma^1(-; \L_{w_1}) 
\ar[r]^-{d} & \hdots \ar[r]^-{d} 
& \Gamma^{n-1}(-; \L_{w_1})
}
\Big)\;.
$$
We first calculate the sheaf cohomology of $Z$ and $\L^{\delta}_{w_1}$ via the Mayer-Vietoris sequence. 
Let ${\cal N}$ be an open tubular neighborhood of the equator in $S^{n}$ and let $W$ be the complement of the closure $\overline{(1-\epsilon){\cal N}}$. In the quotient, this gives a cover $\{U,V\}$ of $\RR P^n$ with $U\simeq B^n$ and $V\simeq \RR P^{n-1}$, where $B^n$ is the $n$-dimensional ball. The intersection $U\cap V\simeq S^{n-1}$ and we have the Mayer-Vietoris sequence with local coefficients
$$ 
\xymatrix{
\ar[r] \hdots  &H^k(\RR P^n; Z)\ar[r] & H^k(B^n;Z)\oplus H^k(\RR P^{n-1};Z)\ar[r] & 
H^{k}(S^{n-1};Z)\ar[r] & \hdots \;.
}
$$ 
Since the restriction of the tautological bundle to $B^n$ and $S^{n-1}$ trivializes, and the restriction to $\RR P^{n-1}$ is the tautological bundle over $\RR P^{n-1}$, the sequence reduces to 
$$ 
\xymatrix{
\ar[r] \hdots  &H^k(\RR P^n; Z)\ar[r] & H^k(\RR P^{n-1};Z)\ar[r] & 
H^{k}(S^{n-1};\ZZ)\ar[r] & \hdots \;.
}
$$ 
Thus, for $n>1$ and $1\leq k\neq n,n-1$, we have an isomorphism
$$H^k(\RR P^n; Z)\cong H^k(\RR P^{n-1};Z)\;.$$
We also have the sequence 
\(\label{rpn seq 2}
\xymatrix{
0\ar[r] & H^{n-1}(\RR P^n;Z)\ar[r] &
H^{n-1}(\RR P^{n-1}; Z)\ar[r] &
 \ZZ\ar[r] & H^n(\RR P^n; Z)\ar[r] & 0\;.
}
\)
We have already shown that $H^1(S^1;Z)\cong \ZZ/2$. The sequence \eqref{rpn seq 2} implies that $
H^{1}(\RR P^2; Z)\cong \ZZ/2$ and $H^2(\RR P^2;Z)\cong \ZZ$.
We claim that for even $n$
$$
H^k(\RR P^n;Z)=\left\{
\begin{array}{cc}
\ZZ , & k=n
\\
\ZZ/2, & 0<k<n \ \ {\rm odd}
\\
0, & {\rm otherwise},
\end{array}
\right.
$$
while for odd $n$
$$
H^k(\RR P^n;Z)=\left\{
\begin{array}{cc}
\ZZ/2, & 0<k\leq n \ \ {\rm odd}
\\
0, & {\rm otherwise}.
\end{array}
\right.
$$
Note the shifts in degrees compared to integral coefficients. 
To prove the claim, we proceed by induction on the dimension $n$. The only nontrivial part of the induction is to show that $H^n(\RR P^n;Z)\cong \ZZ/2$ when $n$ is odd. To prove this, we use the sequence \eqref{rpn seq 2}. By the induction hypothesis, this reduces to the exact sequence
\(\label{ses rpn zs}
\xymatrix{
\ZZ\ar[r] & \ZZ \ar[r] &
H^n(\RR P^n;Z)\ar[r] & 0\;,
}
\)
and we need to identify the map between the copies of $\ZZ$. By the induction by hypothesis, the generator of $H^{n-1}(\RR P^{n-1};Z)$ maps to the generator of $H^{n-1}(S^{n-1};\ZZ)$ under the quotient $q:S^{n-1}\to \RR P^{n-1}$. Consider the commutative diagram
\(\label{rpn mob}
\xymatrix{
S^{n-1}\coprod S^{n-1}\ar[r]&  {\cal N}\simeq  S^{n-1}\ar[d]
\\
S^{n-1} \; \ar@{^{(}->}[r]\ar[u]^{{\rm id}\times -{\rm id}} & M \simeq  \RR P^{n-1}
\;,
}
\)
where ${\cal N}$ is a tubular neighborhood of $S^{n-1}$ in $S^n$ and $-{\rm id}$ the antipodal map. In ordinary cohomology the left map in \eqref{rpn mob} induces the map $(x,y)\mapsto x-y$. However, with $\L$ twisted coefficients, the map is modified by a local trivialization of $\L$ and the resulting map sends $(x,y)\mapsto x+y$. The top map induces the map which sends the generator $x\in H^{n-1}(S^{n-1};\ZZ)$ to $(x,x)\in H^{n-1}(S^{n-1};\ZZ)\oplus H^{n-1}(S^{n-1};\ZZ)$. Thus, the commutativity of the diagram implies that the restriction must send the generator of $H^{n-1}(\RR P^{n-1};Z)$ to twice the generator of $H^{n-1}(S^{n-1};\ZZ)$. Hence the map $\ZZ\to \ZZ$ in \eqref{ses rpn zs} is the $\times 2$ map and $H^n(\RR P^n;Z)\cong \ZZ/2$. This proves the claim.

Similar computations hold for coefficients in $\L^{\delta}$, where the torsion is killed. Using the diagonal sequence
$$
\xymatrix{
H^{k-1}(M;Z)\ar[r] & \Omega^{k-1}(M)/{\rm im}(d)
\ar[r] & \widehat{H}^k(M;\nabla)\ar[r] & H^{k}(M;Z)
\ar[r] & 0
}\;,
$$
in the differential cohomology diamond diagram, we compute
$$
\widehat{H}^k(\RR P^n;w_1)=\left\{
\begin{array}{cc}
\ZZ\oplus \Omega^{n-1}(\RR P^n;\L)/{\rm im(d)} & k=n
\\
\ZZ/2\oplus \Omega^{k-1}(\RR P^n;\L)/{\rm im(d)} & 0<k<n \ \ {\rm odd}
\\
 \Omega^{k-1}(\RR P^n;\L)/{\rm im(d)}  & {\rm otherwise}
\end{array}
\right.
$$
for $n$ even and 
$$
\widehat{H}^k(\RR P^n;w_1)=\left\{
\begin{array}{cc}
\ZZ/2\oplus \Omega^{k-1}(\RR P^n;\L)/{\rm im(d)} & 0<k\leq n \ \ {\rm odd}
\\
 \Omega^{k-1}(\RR P^n;\L)/{\rm im(d)}  & {\rm otherwise}
\end{array}
\right.
$$
for $n$ odd.
\end{example}

%%%%%%%%%%%%%%%

\end{document}